\documentclass[12pt]{article}

\title{A Survey on Coefficients of Cyclotomic Polynomials}
\author{Carlo Sanna}
\date{Politecnico di Torino\\
Department of Mathematical Sciences\\
Corso Duca degli Abruzzi 24, 10129 Torino, Italy\\[2pt]
\texttt{carlo.sanna.dev@gmail.com}}

\usepackage[left=1.15in, right=1.15in, top=.72in, bottom=.72in]{geometry}
\usepackage{amsmath}
\usepackage{amssymb}
\usepackage{amsthm}
\usepackage{bm}
\usepackage{xcolor}
\usepackage[colorlinks=true, pdfencoding=auto]{hyperref}
\usepackage{bookmark}
\usepackage{enumitem}
\setenumerate{label=(\roman*), topsep=0.5em, itemsep=0.5em, leftmargin=3.0em}
\usepackage{graphicx}

\newtheorem{theorem}{Theorem}[section]
\newtheorem{corollary}{Corollary}[section]
\newtheorem{lemma}[theorem]{Lemma}
\newtheorem{conjecture}{Conjecture}[section]

\newlength{\bibitemsep}
\newlength{\bibparskip}
\let\oldthebibliography\thebibliography
\renewcommand\thebibliography[1]{%
  \oldthebibliography{#1}%
  \small
  \setlength{\parskip}{\bibitemsep}%
  \setlength{\itemsep}{\bibparskip}%
}

\DeclareMathOperator*{\rad}{rad}

\begin{document}

\pagenumbering{roman}
\maketitle

\begin{abstract}
Cyclotomic polynomials play an important role in several areas of mathematics and their study has a very long history, which goes back at least to Gauss~(1801).
In~particular, the properties of their coefficients have been intensively studied by several authors, and in the last 10 years there has been a burst of activity in this field of research.
This concise survey attempts to collect the main results regarding the coefficients of the cyclotomic polynomials and to provide all the relevant references to their proofs.
\end{abstract}

\newpage

\tableofcontents

\newpage
\pagenumbering{arabic}

\section{Introduction}

Cyclotomic polynomials play an important role in several areas of mathematics and their study has a very long history, which goes back at least to Gauss~(1801)~\cite{Gauss1801}.

For instance, cyclotomic polynomials appear in: 
the solution of the problem of which regular $n$-gons are constructible with straightedge and compass (Gauss--Wantzel theorem~\cite[p.~46]{MR1483895});
elementary proofs of the existence of infinitely many prime numbers equal to $1$, respectively $-1$, modulo $n$, which is a special case of Dirichlet's theorem on primes in arithmetic progressions~\cite[Sections~48--50]{MR0043111};
Witt's proof~\cite{MR3069571} of Wedderburn's little theorem that every finite domain is a field~\cite[Section~13]{MR1838439};
the \mbox{``cyclotomic criterion''} in the study of primitive divisors of Lucas and Lehmer sequences~\cite{MR1863855};
and lattice-based cryptography~\cite{Mukherjee2016, Park2020112585}.

In~particular, the coefficients of cyclotomic polynomials have been intensively studied by several authors, and in the last 10 years there has been a burst of activity in this field of research.
This concise survey attempts to collect the main results regarding the coefficients of the cyclotomic polynomials and to provide all the relevant references to their proofs.
Previous surveys on this topic were given by Lenstra~(1979)~\cite{MR541398}, Vaughan~(1989)~\cite{MR1203320}, and Thangadurai~(2000)~\cite{MR1802391}.

\paragraph*{Acknowledgments}
The author is grateful to Tsit-Yuen~Lam, Pieter~Moree, and Carl~Pomerance, for providing useful suggestions that improved the quality of this survey.
The author is a member of GNSAGA of INdAM and of CrypTO, the group of Cryptography and Number~Theory of Politecnico di Torino.


\subsection{Definitions and basic facts}

Let $n$ be a positive integer.
The \emph{$n$th cyclotomic polynomial} $\Phi_n(X)$ is defined as the monic polynomial whose roots are the $n$th primitive roots of unity, that is,
\begin{equation}\label{equ:Phi_n-definition}
\Phi_n(X) := \prod_{\substack{1 \,\leq\, k \,\leq\, n \\ \gcd(n,\, k) \,=\, 1}} \left(X - e^{2\pi\mathbf{i}k / n}\right) .
\end{equation}
The word ``cyclotomic'' comes from the ancient Greek words ``cyclo'' (circle) and \mbox{``tomos''} (cutting), and refers to how the $n$th roots of unit divide the circle into equal parts.  
Note that, incidentally, the Greek letter $\Phi$ looks a bit like a cut circle.
The degree of $\Phi_n(X)$ is equal to $\varphi(n)$, where $\varphi$ is the Euler totient function.
Despite its definition in terms of complex numbers, it can be proved that $\Phi_n(X)$ has integer coefficients.
Furthermore, $\Phi_n(X)$ is irreducible over $\mathbb{Q}$ and, consequently, it is the minimal polynomial of any primitive $n$th root of unity.
The irreducibility of $\Phi_n(X)$ for $n$ prime was first proved by Gauss~(1801)~\cite{Gauss1801}, and the irreducibility of $\Phi_n(X)$ in general was first proved by Kronecker~(1854)~\cite{Kronecker1854}.
Weintraub~(2013)~\cite{MR3063118} presented proofs of the irreducibility of $\Phi_n(X)$ due to Gauss, Kronecker, Sch\"onemann, and Eisenstein, for $n$ prime, and Dedekind, Landau, and Schur,\footnote{Perhaps curiously, Schur's proof of the irreducibility of $\Phi_n(X)$ was set to rhymes~\cite[pp.~38--41]{Cremer}.} for every $n$.

From~\eqref{equ:Phi_n-definition} it follows easily that
\begin{equation}\label{equ:Xn-1}
X^n - 1 = \prod_{d \,\mid\, n} \Phi_d(X) ,
\end{equation}
which in turn, by the M\"obius inversion formula, yields that
\begin{equation}\label{equ:Phi_n-Mobius}
\Phi_n(X) = \prod_{d \,\mid\, n} \left(X^{n / d} - 1\right)^{\mu(d)} = \prod_{d \,\mid\, n} \left(X^d - 1\right)^{\mu(n / d)} ,
\end{equation}
where $\mu$ is the M\"obius function.
In~particular, we have that
\begin{equation}\label{equ:Phi_p}
\Phi_p(X) = \frac{X^p - 1}{X - 1} = X^{p - 1} + \cdots + X + 1 ,
\end{equation}
for every prime number $p$.

The next lemma collects some important elementary identities, which can be proved either using~\eqref{equ:Phi_n-Mobius} or checking that both sides have the same zeros~\cite{MR541398, MR1802391}.

\begin{lemma}\label{lem:basic}
For every positive integer $n$ and every prime number $p$, we have that:
\begin{enumerate}
\item\label{lem:basic:1} $\Phi_{pn}(X) = \Phi_n(X^p)$ if $p \mid n$;
\item\label{lem:basic:2} $\Phi_{pn}(X) = \Phi_n(X^p) / \Phi_n(X)$ if $p \nmid n$;
\item\label{lem:basic:3} $\Phi_{2n}(X) = (-1)^{\varphi(n)} \Phi_n(-X)$ if $2 \nmid n$;
\item\label{lem:basic:4} $\Phi_n(X) = \Phi_{\rad(n)}(X^{n / \!\rad(n)})$, where $\rad(n)$ is the product of the primes dividing $n$;
\item\label{lem:basic:5} $\Phi_n(1/X) = X^{-\varphi(n)} \Phi_n(X)$ if $n > 1$.
\end{enumerate}
\end{lemma}

Starting from $\Phi_1(X) = X - 1$ and using Lemma~\ref{lem:basic}'s \ref{lem:basic:1} and \ref{lem:basic:2}, one can inductively compute the cyclotomic polynomials.
The first ten cyclotomic polynomials are:
\begin{equation*}
\begin{aligned}
\Phi_{1}(X) &= X - 1	&	\Phi_{6}(X) &= X^{2} - X + 1	\\
\Phi_{2}(X) &= X + 1	&	\Phi_{7}(X) &= X^{6} + X^{5} + X^{4} + X^{3} + X^{2} + X + 1	\\
\Phi_{3}(X) &= X^{2} + X + 1	&	\Phi_{8}(X) &= X^{4} + 1	\\
\Phi_{4}(X) &= X^{2} + 1	&	\Phi_{9}(X) &= X^{6} + X^{3} + 1	\\
\Phi_{5}(X) &= X^{4} + X^{3} + X^{2} + X + 1	&	\Phi_{10}(X) &= X^{4} - X^{3} + X^{2} - X + 1
\end{aligned}
\end{equation*}
A natural observation is that the coefficients of the cyclotomic polynomials are very small, and one could be even tempted to conjecture that they are always in $\{-1, 0, +1\}$.
The first counterexample to this conjecture occurs for $n = 105$, since we have
\begin{align*}
\Phi_{105}&(X) = X^{48} + X^{47} + X^{46} - X^{43} - X^{42} - \mathbf{2}X^{41} - X^{40} - X^{39} + X^{36} + X^{35} + X^{34} \\
&+ X^{33} + X^{32} + X^{31} - X^{28} - X^{26} - X^{24} - X^{22} - X^{20} + X^{17} + X^{16} + X^{15} \\
&+ X^{14} + X^{13} + X^{12} - X^9 - X^8 - \textbf{2}X^7 - X^6 - X^5 + X^2 + X + 1 .
\end{align*}
It is no coincidence that $105 = 3 \cdot 5 \cdot 7$ is the smallest odd positive integer having three different prime factors.
Indeed, every cyclotomic polynomial $\Phi_n(X)$ such that $n$ has less than three odd prime factors has all its coefficients in $\{-1, 0, +1\}$ (see Section~\ref{sec:binary}).

For every positive integer $n$, let us write
\begin{equation*}
\Phi_n(X) = \sum_{j \,\geq\, 0} a_n(j) X^j , \quad a_n(j) \in \mathbb{Z} ,
\end{equation*}
so that $a_n(j)$ is the coefficient of $X^j$ in $\Phi_n(X)$.
(Note that $a_n(j) = 0$ for $j \notin [0, \varphi(n)]$.)
The peculiarity of the smallness of the coefficients of the cyclotomic polynomials was very well explained by D.~H.~Lehmer~(1966)~\cite{MR197441}, who wrote: ``The smallness of $|a_n(j)|$ would appear to be one of the fundamental conspiracies of the primitive $n$th roots of unity.
When one considers that $a_n(j)$ is a sum of $\binom{\varphi(n)}{j}$ unit vectors (for example $73629072$ in the case of $n = 105$, $j = 7$) one realizes the extent of the cancellation that takes place.''

In light of Lemma~\ref{lem:basic}'s \ref{lem:basic:3} and \ref{lem:basic:4}, for the purpose of studying the coefficients of $\Phi_n(X)$ it suffices to consider only odd squarefree integers $n$.
A squarefree positive integer $n$, or a cyclotomic polynomial $\Phi_n(X)$, is \emph{binary}, \emph{ternary}, \ldots if the number of prime factors of $n$ is $2$, $3$, \ldots, respectively.
The \emph{order} of $\Phi_n(X)$ is the number of prime factors of $n$.
From Lemma~\ref{lem:basic}'s \ref{lem:basic:5} we have that for every integer $n > 1$ the cyclotomic polynomial $\Phi_n(X)$ is \emph{palindromic}, that is,
\begin{equation*}
a_n(\varphi(n) - j) = a_n(j), \quad \text{ for } j = 0,\dots,\varphi(n) .
\end{equation*}

We conclude this section by defining the main quantities that have been considered in the study of the coefficients of the cyclotomic polynomials. 
First, we have $A(n)$, $A^+(n)$, and $A^-(n)$, which are defined as follows
\begin{equation*}
A(n) := \max_{j \,\geq\, 0} |a_k(j)|, \quad A^+(n) := \max_{j \,\geq\, 0} a_k(j), \quad A^-(n) := \min_{j \,\geq\, 0} a_k(j) .
\end{equation*}
In general, the \emph{height} of a polynomial $P \in \mathbb{C}[X]$ is defined as the maximum of the absolute values of the coefficients of $P$, and $P$ is \emph{flat} if its height is not exceeding $1$.
Thus, $A(n)$ is the height of $\Phi_n(X)$.
We also let $\mathcal{A}(n) := \{a_n(j) : 0 \leq j \leq \varphi(n)\}$ be the set of coefficients of $\Phi_n(X)$.
Moreover, we let $\theta(n) := \#\{j \geq 0 : a_n(j) \neq 0\}$ be the number of nonzero coefficients of $\Phi_n(X)$.
The \emph{maximum gap} of a nonzero polynomial $P(X) = \sum_{i=1}^k c_k X^{e_k} \in \mathbb{C}[X]$, where $c_1, \dots, c_k \in \mathbb{C}^*$ and $e_1 < \cdots < e_k$, is defined as $G(P) := \max\{e_{j+1} - e_j : j < k\}$.
We let $G(n) := G(\Phi_n)$ denote the maximum gap of $\Phi_n(X)$.
Note that by~\eqref{equ:Phi_p} we have that $A(p) = A^+(p) = A^-(p) = 1$, $\mathcal{A}(p) = \{1\}$, $\theta(p) = p$, and $G(p) = 1$, for every prime number $p$.
Thus, the first interesting case in the study of these quantities is the one of binary cyclotomic polynomials.

\section{Binary cyclotomic polynomials}\label{sec:binary}

The understanding of the coefficients of binary cyclotomic polynomials is quite complete.
Let $p$ and $q$ be distinct prime numbers.
From~\eqref{equ:Phi_n-Mobius} it follows that
\begin{equation*}
\Phi_{pq}(X) = \frac{(X^{pq} - 1)(X - 1)}{(X^p - 1)(X^q - 1)} .
\end{equation*}
Migotti~(1883)~\cite{zbMATH02705050} and Bang~(1895)~\cite{zbMATH02678148} proved that the coefficients of every binary cyclotomic polynomial belong to $\{+1, -1, 0\}$.
Beiter (1964)~\cite{MR1532827} gave a first criterion to enstablish if $a_{pq}(j)$ is equal to $+1$, $-1$, or $0$.
This criterion is a bit difficult to apply, but she used it to compute the midterm coefficient of $\Phi_{pq}(X)$.
Using a different method, Habermehl, Richardson, and Szwajkos~(1964)~\cite{MR1571408} determined the coefficients of $\Phi_{3p}(X)$, for $p > 3$.
Carlitz~(1966)~\cite{MR202655} gave a formula for the number of nonzero coefficients of $\Phi_{pq}(X)$, and Lenstra~(1979)~\cite{MR541398} proved an expansion for $\Phi_{pq}(X)$, which was then rediscovered by Lam and Leung~(1996)~\cite{MR1404079}, that leads to an explicit determination of $a_{pq}(j)$.
Moree~(2014)~\cite{MR3295662} generalized this formula to binary inclusion-exclusion polynomials (see~Section~\ref{sec:inclusion-exclusion}), and he also showed a connection with numerical semigroups.

The following theorem gives a precise description of the coefficients of binary cyclotomic polynomials~\cite{MR1404079, MR541398, MR3295662, MR1802391}.

\begin{theorem}\label{thm:binary}
Let $p < q$ be distinct prime numbers, and let $\overline{p}$ and $\overline{q}$ be the unique positive integers such that $pq + 1 = p\overline{p} + q\overline{q}$.
(Equivalently, $\overline{p}$ is the inverse of $p$ modulo $q$ and $\overline{q}$ is the inverse of $q$ modulo $p$.)
We have that:
\begin{enumerate}
\item It holds
\begin{equation*}
\Phi_{pq}(X) = \sum_{i \,=\, 0}^{\overline{p} - 1} X^{pi} \sum_{j \,=\, 0}^{\overline{q} - 1} X^{qj} - X^{-pq} \sum_{i \,=\, \overline{p}}^{q - 1} X^{pi} \sum_{j \,=\, \overline{q}}^{p - 1} X^{qj} .
\end{equation*}

\item\label{thm:binary:2} For every nonnegative integer $j < pq$, we have that either $j = px + qy$ or $j = px + qy - pq$ with $x < q$ the unique nonnegative integer such that $px \equiv j \pmod q$ and $y < p$ the unique nonnegative integer such that $qy \equiv j \pmod p$; and it holds
\begin{equation*}
a_{pq}(j) = 
    \begin{cases}
    +1 & \text{ if $j = px + qy$ with $0 \leq x < \overline{p}$, $0 \leq y < \overline{q}$}; \\
    -1 & \text{ if $j = px + qy - pq$ with $\overline{p} \leq x < q$, $\overline{q} \leq y < p$}; \\
     0 & \text{ otherwise}.
    \end{cases}
\end{equation*}

\item\label{thm:binary:3} The number of positive coefficients of $\Phi_{pq}(X)$ is equal to $\overline{p}\,\overline{q}$, the number of negative coefficients is equal to $\overline{p}\,\overline{q} - 1$, and (thus) the number of nonzero coefficients of $\Phi_{pq}(X)$ is equal to $2\overline{p}\,\overline{q} - 1$.

\item The nonzero coefficients of $\Phi_{pq}(X)$ alternates between $+1$ and $-1$.

\item The midterm coefficient of $\Phi_{pq}(X)$ satisfies $a_{pq}(\varphi(pq)/2) = (-1)^{\overline{p} \,- 1}$.
\end{enumerate}
\end{theorem}

Moree~\cite{MR3295662} gave a nice way to illustrate Theorem~\ref{thm:binary}'s \ref{thm:binary:2} by using what he called an \emph{LLL-diagram} (for Lenstra, Lam, and Leung).
This is a $p \times q$ matrix constructed as follows.
Start with $0$ in the bottom-left entry, add $p$ for every move to the right, add $q$ for every move upward, and reduce all entries modulo $pq$.
The numbers in the bottom-left $\overline{p} \times \overline{q}$ submatrix are the exponents of the positive terms of $\Phi_{pq}(X)$, and the numbers in the top-right $(p - \overline{p}) \times (q - \overline{q})$ submatrix are the exponents of the negative terms of $\Phi_{pq}(X)$.
For example, the LLL-diagram for the binary cyclotomic polynomial
\begin{align*}
\Phi_{5 \;\!\cdot\;\! 7}(X) = X^{24} &- X^{23} + X^{19} - X^{18} + X^{17} - X^{16} + X^{14} - X^{13} \\
& + X^{12} - X^{11} + X^{10} - X^8 + X^7 - X^6 + X^5 - X + 1
\end{align*}
is the following
\begin{center}
\begin{tabular}{ccc|cccc}
\cline{4-7}
\rule{0pt}{0.8\normalbaselineskip}$28$ & $33$ &  $3$ &  $8$ & $13$ & $18$ & \multicolumn{1}{c|}{$23$} \\
$21$ & $26$ & $31$ &  $1$ &  $6$ & $11$ & \multicolumn{1}{c|}{$16$} \\
\hline
\multicolumn{1}{|c}{$14$}\rule{0pt}{0.8\normalbaselineskip} & $19$ & $24$ & $29$ & $34$ &  $4$ &  $9$ \\
\multicolumn{1}{|c}{ $7$} & $12$ & $17$ & $22$ & $27$ & $32$ &  $2$ \\
\multicolumn{1}{|c}{ $0$} &  $5$ & $10$ & $15$ & $20$ & $25$ & $30$ \\
\cline{1-3}
\end{tabular}
\end{center}

By Theorem~\ref{thm:binary}'s \ref{thm:binary:3}, for every binary number $n = pq$, with $p < q$ primes, the number of nonzero coefficients of $\Phi_n(X)$ is $\theta_n = 2\overline{p}\,\overline{q} - 1$.
From $pq + 1 = p\overline{p} + q\overline{q}$ it follows in an elementary way that $\theta_n > n^{1/2}$~\cite[Section~3.1]{MR3101077}.
Lenstra~(1979)~\cite{MR541398} proved that for every $\varepsilon > 0$ there exist infinitely many binary $n = pq$ such that $\theta_n < p^{8/13 + \varepsilon}$.
The proof is based on a result of Hooley~(1973)~\cite{MR354579} that says that for every $\varepsilon > 0$ there exist inﬁnitely many prime number $p$ such that $P(p - 1) > p^{5 / 8 - \varepsilon}$, where $P(n)$ denotes the largest prime factor of $n$.
Hooley's result has been improved by several authors.
Currently, the best bound is $P(p - 1) > p^{0.677}$, which is due to Baker and Harman~(1998)~\cite{MR1610553}.
This reduces the exponent $8 / 13$ of Lenstra's bound to $1 /(1 + 0.677) = 0.596\dots$.
Using a different method, Bzd\c{e}ga~(2012)~\cite{MR2875347} proved that there are infinitely many binary numbers $n$ such that $\theta_n < n^{1/2 + \varepsilon}$, and also gave upper and lower bounds for the number $H_\varepsilon(x)$ of binary $n \leq x$ such that $\theta_n < n^{1/2 + \varepsilon}$.
Fouvry~(2013)~\cite{MR3101077} proved the following asymptotic formula for $H_\varepsilon(x)$.

\begin{theorem}
For $\varepsilon \in (0, 1/2)$, let
\begin{equation*}
C(\varepsilon) := \frac{2}{1 + 2\varepsilon} \log\!\left(\frac{1 + 2\varepsilon}{1 - 2\varepsilon}\right) .
\end{equation*}
Then for every $\varepsilon_0 > 0$, uniformly for $\varepsilon \in (12/15 + \varepsilon_0, 1/2 - \varepsilon_0)$, we have that
\begin{equation}\label{equ:binary-Fouvry}
H_\varepsilon(x) \sim C(\varepsilon) \, \frac{x^{1/2 + \varepsilon}}{\log x} 
\end{equation}
as $x \to +\infty$.
\end{theorem}

Furthermore, Fouvry~\cite{MR3101077} provided an upper bound and a lower bound for $H_\varepsilon(x)$ of the same order of~\eqref{equ:binary-Fouvry}, and showed that the Elliott--Halberstam Conjecture implies that~\eqref{equ:binary-Fouvry} holds in the range $\varepsilon \in (\varepsilon_0, 1/2 - \varepsilon_0)$.

Hong, Lee, Lee, and Park~(2012)~\cite{MR2944756} determined the maximum gap of binary cyclotomic polynomials, and Moree~(2014)~\cite{MR3295662} gave another proof of the result using numerical semigroups.
Yet another short proof was given by Kaplan~(2016)~\cite[End of Section~2.1]{MR3459568}.
Furthermore, Camburu, Ciolan, Luca, Moree, and Shparlinski~(2016)~\cite{MR3459568} determined the number of maximum gaps of $\Phi_{pq}(X)$, and the existence of particular gaps in the case in which $q \equiv \pm 1 \pmod p$.
The following theorem collects these results~\cite{MR3459568, MR2944756, MR3295662}.

\begin{theorem}\label{thm:binary-gap}
Let $p < q$ be prime numbers.
Then:
\begin{enumerate}
\item $G(pq) = p - 1$.
\item The number of maximum gaps of $\Phi_{pq}(X)$ is equal to $2 \lfloor q / p\rfloor$.
\item $\Phi_{pq}(X)$ contains the sequence of consecutive coefficients
\begin{equation*}
\pm 1, \underbrace{0, \dots, 0}_{\text{$m$ times}}, \pm 1
\end{equation*}
for all $m \in \{0, \dots, p - 2\}$ if and only if $q \equiv \pm 1 \pmod p$.
\end{enumerate}
\end{theorem}

Cafure and Cesaratto~(2021)~\cite{XX0001} considered the coefficients of $\Phi_{pq}(X)$ as a word over the ternary alphabet $\{+1, -1, 0\}$, and provided an algorithm that, given as input $p < q$ and the quotient and remainder of the division of $q$ by $p$, computes $\Phi_{pq}(X)$ performing $O(pq)$ simple operations on words.
Chu~(2021)~\cite{MR4217760} proved that the exponents of the positive, respectively negative, terms of $\Phi_{pq}(X)$ are in arithmetic progression if and only if $q \equiv 1 \pmod p$, respectively $q \equiv -1 \pmod p$.

\section{Ternary cyclotomic polynomials}

Ternary cyclotomic polynomials are the simplest ones for which the behavior of the coefficients is not completely understood.
Kaplan~(2007)~\cite[Lemma~1]{MR2351667} proved the following lemma, which provides a formula for the coefficients of ternary cyclotomic polynomials.
This is known as \emph{Kaplan's lemma} and has been used to prove several results on ternary cyclotomic polynomials~\cite{MR3208875, MR2568054, MR2544145, MR2905233, MR2718826, MR2978846, MR3248694, MR3614561, MR3606634, MR4261643, MR3432726}.

\begin{lemma}[Kaplan's lemma]\label{lem:Kaplan-lemma}
Let $p < q < r$ be odd prime numbers and let $j \geq 0$ be an integer.
For every integer $i \in {[0, pq)}$, put
\begin{equation*}
b_i := 
    \begin{cases}
    a_{pq}(i) & \text{ if $ri \leq j$}; \\
    0 & \text{ otherwise}.
    \end{cases}
\end{equation*}
Then we have
\begin{equation*}
a_{pqr}(j) = \sum_{m \,=\, 0}^{p - 1} (b_{f(m)} - b_{f(m + q)}) ,
\end{equation*}
where $f(m)$ is the unique integer such that $f(m) \equiv r^{-1}(j - m) \pmod {pq}$, $0 \leq f(m) < pq$.
\end{lemma}

Lemma~\ref{lem:Kaplan-lemma} reduces the computation of $a_{pqr}(j)$ to that of $a_{pq}(i)$, which in turn is provided by Theorem~\ref{thm:binary}'s~\ref{thm:binary:2}.
Note that in order to compute $a_{pqr}(j)$ using Lemma~\ref{lem:Kaplan-lemma} and Theorem~\ref{thm:binary}'s~\ref{thm:binary:2} it is not necessary to compute $f(m)$ and $f(m + q)$, but it suffices to compute $[f(m)]_p$, $[f(m)]_q$, $[f(m + q)]_p$, and $[f(m + q)]_q$, which can be easier, where $[k]_p$ and $[k]_q$ are the unique nonnegative integers $x < q$ and $y < p$ such that $px \equiv k \pmod q$ and $qy \equiv k \pmod p$, for every integer $k$.
Actually, since $[f(m)]_p = [f(m + q)]_p$, it suffices to compute $[f(m)]_p$, $[f(m)]_q$, and $[f(m + q)]_q$.

For the rest of this section, let $p < q < r$ be odd prime numbers and let $n = pqr$ be a ternary integer.
The next subsections describe the main themes of research on ternary cyclotomic polynomials.

\subsection{Bounds on the height and Beiter's conjecture}

Upper bounds for the height of ternary cyclotomic polynomials have been studied by many authors.
Bang~(1895)~\cite{zbMATH02678148} proved that $A(pqr) \leq p - 1$.
Beiter~(1968)~\cite{MR227085} made the following conjecture, which is known as \emph{Beiter's conjecture}.

\begin{conjecture}[Beiter's conjecture]\label{con:Beiter-conjecture}
$A(pqr) \leq \tfrac1{2}(p + 1)$ for all odd primes $p < q < r$.
\end{conjecture}

Beiter~(1968)~\cite{MR227085} proved her conjecture in the case in which $q \equiv \pm 1 \pmod p$ or $r \equiv \pm 1 \pmod p$.
As~a consequence, Beiter's conjecture is true for $p = 3$.
Also, Bloom~(1968)~\cite{MR227086} showed that Beiter's conjecture is true for $p = 5$.
Beiter~(1971)~\cite{MR280432} improved Bang's bound to $A(pqr) \leq p - \lfloor (p + 1) / 4 \rfloor$.
M\"oller~(1971)~\cite{MR274383} proved that for every odd prime number $p$ there exists a ternary cyclotomic polynomial $\Phi_{pqr}(X)$, with $p < q < r$, having a coefficient equal to $(p + 1)/2$.
This shows that Beiter's conjecture, if true, is the best possible.
Bachman~(2003)~\cite{MR1971249} proved an upper bound for $A^+(pqr)$ and a lower bound for $A^-(pqr)$ in terms of $p$ and the inverses of $q$ and $r$ modulo $p$.
As corollaries, he deduced that: Beiter's conjecture is true if $q$ or $r$ is equal to $\pm 1, \pm 2$ modulo $p$; we have $A(pqr) \leq p - \lceil p / 4 \rceil$; and $A^+(pqr) - A^-(pqr) \leq p$, in~particular either $A^+(pqr) \leq (p - 1) / 2$ or $A^-(pqr) \geq -(p - 1) / 2$.
Note that the first two corollaries improve the previous results of Beiter~\cite{MR227085, MR280432}.
Regarding the third, Bachman~(2004)~\cite{MR2053964} also proved that for every odd prime number $p$ there exist infinitely many ternary cyclotomic polynomials $\Phi_{pqr}(X)$, with $p < q < r$, such that $\mathcal{A}(pqr) = [-(p - 1)/2, (p + 1)/2] \cap \mathbb{Z}$, and similarly for the interval $[-(p + 1)/2, (p - 1)/2]$.
Leher~(2007)~\cite[p.~70]{Leher2007} found a counterexample to Beiter's conjecture, that is, $A(17 \cdot 29 \cdot 41) = 10 > (17 + 1) / 2$.
Let $M(p) := \max_{p < q < r} A(p)$.
For every odd prime $p$, Gallot and Moree~(2009)~\cite{MR2544145} defined an effectively computable set of natural numbers $\mathcal{B}(p)$ such that if $\mathcal{B}(p)$ is nonempty then 
\begin{equation*}
M(p) \geq p - \min(\mathcal{B}(p)) > (p + 1)/2, 
\end{equation*}
and so Beiter's conjecture is false for $p$.
Then, for $p \geq 11$ they showed that $\mathcal{B}(p)$ is nonempty and $\max(\mathcal{B}(p)) = (p - 3)/2$.
Moreover, for every $\varepsilon > 0$, they proved that
\begin{equation}\label{equ:Mp-bounds}
\left(\tfrac{2}{3} - \varepsilon\right)p \leq M(p) \leq \tfrac{3}{4} p ,
\end{equation}
for all sufficiently large $p$.
In light of these results, they formulated the following:

\begin{conjecture}[Corrected Beiter's conjecture]
$M(p) \leq \tfrac{2}{3}p$ for every prime $p$.
\end{conjecture}

Zhao and Zhang~(2010)~\cite{MR2660889} gave a sufficient condition for the Corrected Beiter conjecture and proved it when $p = 7$.
(Note that for $p = 7$ the Corrected Beiter Conjecture is equivalent to the original Beiter Conjecture.)
Moree and Ro\c{s}u~(2012)~\cite{MR2978846} showed that for each odd integer $\ell \geq 1$ there exist infinitely many odd primes $p < q < r$ such that 
\begin{equation*}
\mathcal{A}(pqr) = [-(p-\ell-2)/2, (p+\ell+2)/2] \cap \mathbb{Z} .
\end{equation*}
This provides a family of cyclotomic polynomials that contradict the Beiter conjecture and have the largest coefficient range possible.
Bzd\c{e}ga~(2010)~\cite{MR2660553} improved Bachman's bounds~\cite{MR1971249} by giving the following theorem.

\begin{theorem}\label{thm:Bzd-bounds}
Let $p < q < r$ be odd primes and let $q^\prime$ and $r^\prime$ be the inverses of $q$ and $r$ modulo $p$, respectively.
Then
\begin{equation*}
A^+(pqr) \leq \min\{2\alpha + \beta, p - \beta\}, \quad -A^-(pqr) \leq \min\{p + 2\alpha - \beta, \beta\} ,
\end{equation*}
\begin{equation*}
A(pqr) \leq \min\{2\alpha + \beta^*, p - \beta^*\} ,
\end{equation*}
where $\alpha := \min\{q^\prime, r^\prime, p - q^\prime, p - r^\prime\}$, $\beta$ is the inverse of $\alpha qr$ modulo $p$, and $\beta^* := \min\{\beta, p - \beta\}$.
\end{theorem}
As an application of Theorem~\ref{thm:Bzd-bounds}, Bzd\c{e}ga proved a density result showing that the Corrected Beiter conjecture holds for at least $25/27 + O(1/p)$ of all the ternary cyclotomic polynomials with the smallest prime factor dividing their order equal to $p$. 
He also proved that for these polynomials the average value of $A(pqr)$ does not exceed $(p + 1)/2$.
Moreover, for every prime $p \geq 13$, he provided some new classes of ternary cyclotomic polynomials $\Phi_{pqr}(X)$ for which the set of coefficients is very small.
Luca, Moree, Osburn, Saad~Eddin, and Sedunova~(2019)~\cite{MR3922607}, using Theorem~\ref{thm:Bzd-bounds} and some analytic estimates for constrained ternary integers that they developed, showed that the relative density of ternary integers for which the correct Sister Beiter conjecture holds true is at least $25/27$.

Gallot, Moree, and Wilms~(2011)~\cite{MR2905233} initiated the study of
\begin{equation*}
M(p, q) := \max\{A(pqr) : r > q\} .
\end{equation*}
They remarked that $M(p, q)$ can be effectively computed for any given odd primes $p < q$.
For $p = 3, 5, 7, 11, 13, 19$, they proved that the set $\mathcal{Q}_p$ of primes $q$ with $M(p, q) = M(p)$ has a subset of positive density, which they determined, and they also conjectured the value of the natural density of $\mathcal{Q}_p$.
Moreover, they computed or bound $M(p, q)$ for $p$ and $q$ satisfying certain conditions, and they posed several problems regarding $M(p, q)$~\cite[Section~11]{MR2905233}.
Cobeli, Gallot, Moree, and Zaharescu~(2013)~\cite{MR3124808}, using techniques from the study of the distribution of modular inverses, in particular bounds on Kloosterman sums, improved the lower bound in~\eqref{equ:Mp-bounds} to
\begin{equation*}
M(p) > \tfrac{2}{3}p - 3p^{3/4} \log p ,
\end{equation*}
for every prime $p$, and
\begin{equation*}
M(p) > \tfrac{2}{3}p - C p^{1/2} ,
\end{equation*}
for infinitely many primes $p$, where $C > 0$ is a constant.
Moreover, they proved that
\begin{equation*}
\liminf_{x \to +\infty} \frac{\#\{q : p < q \leq x, \; M(p, q) > (p + 1) / 2 \}}{\#\{p : p \leq x\}} \geq \frac{\#\mathcal{B}(p)}{p - 1} .
\end{equation*}
Duda~(2014)~\cite{MR3208875} put $M_{q^\prime}(p) := \max\{M(p, q) : q \equiv q^\prime \pmod p\}$ and proved one of the main conjectures on $M(p, q)$ of Gallot, Moree, and Wilms~\cite[Conjecture~8]{MR2905233}, that is, for all distinct primes $p$ and $q^\prime$ there exists $q_0 \equiv q^\prime \pmod p$ such that for every prime $q \geq q_0$ with $q \equiv q^\prime \pmod p$ we have $M(p, q) = M_{q^\prime}(p)$.
Also, he gave an effective method to compute $M_q(p)$, from which it follows an algorithm to effectively compute $M(p)$, since $M(p) = \max\{M_q(p) : q < p\}$.

Kosyak, Moree, Sofos, and Zhang~(2021)~\cite{MR4196796} conjectured that every positive integer is of the form $A(n)$, for some ternary integer $n$.
They proved this conjecture under a stronger form of Andrica's conjecture on prime gaps, that is, assuming that $p_{n+1} - p_n < \sqrt{p_n} + 1$ holds for every $n \geq 31$, where $p_n$ denotes the $n$th prime number.
Furthermore, they showed that almost all positive integers are of the form $A(n)$ where $n = pqr$ with $p < q < r$ primes is a ternary integer and $\#\mathcal{A}(n) = p + 1$ (which is the maximum possible value for this cardinality).
A nice survey regarding these connections between cyclotomic polynomials and prime gaps was given by Moree~(2021)~\cite{Moree2021}.

\subsection{Flatness}

Recall that a cyclotomic polynomial $\Phi_n(X)$ is \emph{flat} if $A(n) = 1$.
Several families of flat ternary cyclotomic polynomials have been constructed, but a complete classification is still not known.
Beiter~(1978)~\cite{MR514317} characterized the primes $r > q > 3$ such that $\Phi_{3qr}(X)$ is flat.
In~particular, there are infinitely many such primes.
Bachman~(2006)~\cite{MR2201603} proved that if $p \geq 5$, $q \equiv -1 \pmod p$, and $r \equiv 1 \pmod{pq}$ then $\Phi_{pqr}(X)$ is flat.
Note that, for every prime $p \geq 5$, the existence of infinitely many primes $q$ and $r$ satisfying the aforementioned congruences is guaranteed by Dirichlet's theorem on primes in arithmetic progressions.
Flanagan~(2007)~\cite{Flanagan2007} improved Bachman's result by relaxing the congruences to $q \equiv \pm 1 \pmod p$ and $r \equiv \pm 1 \pmod{pq}$.
Kaplan~(2007)~\cite{MR2351667} used Lemma~\ref{lem:Kaplan-lemma} to show that last congruence suffices, that is, the following holds:

\begin{theorem}\label{thm:flat-ternary}
$\Phi_{pqr}(X)$ is flat for all primes $p < q < r$ with $r \equiv \pm 1 \pmod {pq}$.
\end{theorem}

Luca, Moree, Osburn, Saad~Eddin, and Sedunova~(2019)~\cite{MR3922607} proved some asymptotic formulas for ternary integers that, together with Theorem~\ref{thm:flat-ternary}, yield that for every sufficiently large $N > 1$ there are at least $C N / \log N$ ternary integers $n \leq N$ such that $\Phi_n(X)$ is flat, where $C := 1.195\dots$ is an explicit constant.

Ji~(2010)~\cite{MR2718826} considered odd primes $p < q < r$ such that $2r \equiv \pm 1 \pmod {pq}$ and showed that in such a case $\Phi_{pqr}(X)$ is flat if and only if $p = 3$ and $q \equiv 1 \pmod 3$.
For $a \in \{3, 4, 5\}$, Zhang~(2017)~\cite{MR3606634} gave similar characterizations for the odd primes such that $ar \equiv \pm 1 \pmod {pq}$ and $\Phi_{pqr}(X)$ is flat (see also~\cite{MR3432726} for a weaker result for the case $a = 4$).
For $a \in \{6, 7\}$, Zhang~(2020, 2021)~\cite{MR4128150, MR4261643} characterized the odd primes such that $q \equiv \pm 1 \pmod p$, $ar \equiv \pm 1 \pmod {pq}$, and $\Phi_{pqr}(X)$ is flat.
Zhang~(2017)~\cite{MR3645283} also showed that if $p \equiv 1 \pmod w$, $q \equiv 1 \pmod {pw}$, and $r \equiv w \pmod {pq}$, for some integer $w \geq 2$, then $A^+(pqr) = 1$.
(See also the unpublished work of Elder~(2012)~\cite{Elder2012}.)
Furthermore, for $q \not\equiv 1 \pmod p$ and $r \equiv -2 \pmod {pq}$, Zhang~(2014)~\cite{MR3248694} constructed an explicit $j$ such that $a_{pqr}(j) = -2$, so that $\Phi_{pqr}(X)$ is not flat.
Regarding nonflat ternary cyclotomic polynomials with small heights, Zhang~(2017)~\cite{MR3614561} showed that for every prime $p \equiv 1 \pmod 3$ there exist infinitely many $q$ and $r$ such that $A(pqr) = 3$.

\subsection{Jump one property}

Gallot and Moree~(2009)~\cite{MR2568054} proved that neighboring coefficients of ternary cyclotomic polynomials differ by at most one.
They called this property \emph{jump one property}.

\begin{theorem}[Jump one property]
Let $n$ be a ternary integer.
Then
\begin{equation*}
|a_n(j) - a_n(j - 1)| \leq 1
\end{equation*}
for every integer $j \geq 1$.
\end{theorem}

\begin{corollary}
Let $n$ be a ternary integer.
Then $\mathcal{A}(n)$ is a set of consecutive integers.
\end{corollary}

Gallot and Moree used the \emph{jump one property} to give a different proof of Bachman's result~\cite{MR2053964} on ternary polynomials with optimally large set of coefficients.
Their proof of the jump one property makes use of Kaplan's lemma.
Previously, Leher~(2007)~\cite[Theorem~57]{Leher2007} proved the bound $|a_n(j) - a_n(j - 1)| \leq 4$ using methods from the theory of numerical semigroups.
A different proof of the jump one property was given by Bzd\c{e}ga~(2010)~\cite{MR2660553}.
Furthermore, for every ternary integer $n$, Bzd\c{e}ga~(2014)~\cite{MR3206393} gave a characterization of the positive integers $j$ such that $|a_n(j) - a_n(j - 1)| = 1$.
A coefficient $a_n(j)$ is \emph{jumping up}, respectively \emph{jumping down}, if $a_n(j) = a_n(j - 1) + 1$, respectively $a_n(j) = a_n(j - 1) - 1$.
Since cyclotomic polynomials are palindromic, the number of jumping up coefficients is equal to the number of jumping down coefficients.
Let $J_n$ denote such number.
Bzd\c{e}ga~\cite{MR3206393} proved that $J_n > n^{1/3}$ for all ternary integers $n$.
As a corollary, $\theta_n > n^{1/3}$.
Also, he showed that Schinzel~Hypothesis~H implies that for every $\varepsilon > 0$ we have $J_n < 10n^{1/3 + \varepsilon}$ for infinitely many ternary integers $n$.
Camburu, Ciolan, Luca, Moree, and Shparlinski~(2016)~\cite{MR3459568} gave an unconditional proof that $J_n < n^{7/8 + \varepsilon}$ for infinitely many ternary integers $n$.

\section{Higher order cyclotomic polynomials}

There are few specific results regarding cyclotomic polynomials of order greater than three.
Bloom~(1968)~\cite{MR227086} proved that, for odd prime numbers $p < q < r < s$, it holds $A(pqrs) \leq p(p - 1)(pq - 1)$.
Kaplan~(2010)~\cite{MR2684126} constructed the first infinite family of flat cyclotomic polynomials of order four.
Precisely, he proved that $\Phi_{3 \cdot 5 \cdot 31 \cdot s}(X)$ is flat for every prime number $s \equiv -1 \pmod {465}$.
Also, he suggested that all flat cyclotomic polynomials $\Phi_{pqrs}(X)$ satisfy $q \equiv -1 \pmod p$, $r \equiv \pm 1 \pmod {pq}$, and $s \equiv \pm 1 \pmod {pqr}$.
Furthermore, Bzd{\c{e}}ga~(2012)~\cite{MR2890546} proved the upper bounds
\begin{equation*}
A(pqrs) \leq \frac{3}{4}p^3 q, \quad A(pqrst) \leq \frac{135}{512} p^7 q^3 r, \quad A(pqrstu) \leq \frac{18225}{262144} p^{15} q^7 r^3 s ,
\end{equation*}
for all odd prime numbers $p < q < r < s < t < u$.

\section{Height of cyclotomic polynomials}

\subsection{Asymptotic bounds on $A(n)$}

Schur~(1931)\footnote{Unpublished letter to Landau, see~\cite{MR1563307}.} was the first to prove that the coefficients of cyclotomic polynomials can be arbitrarily large, that is, $\sup_{n \geq 1} A(n) = +\infty$.
E.~Lehmer~(1936)~\cite{MR1563307} presented Schur's proof and proved the stronger result that $A(n)$ is unbounded also when $n$ is restricted to ternary integers.
Erd\H{o}s~(1946)~\cite{MR14110} proved that $A(n) > \exp(C(\log n)^{4/3})$ for infinitely many positive integers $n$, for some constant $C > 0$.
His proof rests on a lower bound for the maximum of $|\Phi_n(X)|$ on the unit circle, and the simple consideration that $|\Phi_n(z)| \leq n A(n)$ for every $z \in \mathbb{C}$ with $|z| \leq 1$.
This is essentially the main technique that has then been used to prove lower bounds for $A(n)$~\cite{MR3842902, MR35787, MR91966, MR2106468, MR1084190, MR1225426, MR1830573, MR364141}.
Furthermore, Erd\H{o}s suggested that\footnote{The following formula was printed incorrectly in Erd\H{o}s' paper~\cite{MR14110}, see~\cite{MR32677}.} $A(n) > \exp(n^{C / \log \log n})$ for infinitely many positive integers $n$, for some constant $C > 0$, and claimed that this is the best possible upper bound.
Bateman~(1949)~\cite{MR32677} gave a short proof that, for every $\varepsilon > 0$, it holds $A(n) < \exp(n^{(1+\varepsilon)\log 2/\log \log n})$ for all sufficiently large integers $n$.
Hence, the lower bound suggested by Erd\H{os}, if true, is indeed the best possible.
Then Erd\H{o}s~(1949)~\cite{MR35787} proved that in fact $A(n) > \exp(n^{C/\log \log n})$ for infinitely many positive integers $n$, for some constant $C > 0$, by showing that $\max_{|z| = 1} |\Phi_n(z)| > \exp(n^{C/\log \log n})$ for infinitely many positive integers $n$.
His proof of this last fact is quite involved.
Later, Erd\H{o}s~(1957)~\cite{MR91966} found a simpler proof of the fact that $\max_{x \in (0, 1)} |\Phi_n(x)| > \exp(n^{C/\log \log n})$ for infinitely many positive integers $n$, which again implies the lower bound on $A(n)$.
He conjectured that one can take every positive constant $C < \log 2$, and so Bateman's result is the best possible.
This conjecture was settled by Vaughan~(1974)~\cite{MR364141}, who proved that actually $C = \log 2$ is admissible (see also~\cite{MR3842902} for an alternative proof).

In summary, the maximal order of $A(n)$ is given by the following theorem~\cite{MR32677, MR364141}.
\begin{theorem}[Bateman--Vaughan]
On the one hand, for every $\varepsilon > 0$, we have
\begin{equation*}
A(n) < \exp\!\left(n^{(\log 2 + \varepsilon) / \log \log n}\right)
\end{equation*}
for all sufficiently large positive integers $n$.
On the other hand, we have
\begin{equation*}
A(n) > \exp\!\left(n^{\log 2 / \log \log n}\right)
\end{equation*}
for infinitely many positive integers $n$.
\end{theorem}

Maier~(1990, 1996)~\cite{MR1084190, MR1409383} proved that $n^{f(n)} < A(n) < n^{g(n)}$ for almost all positive integers, where $f$ and $g$ are arbitrary functions such that $f(n) \to 0$ and $g(n) \to +\infty$ as $n \to +\infty$.
Furthermore, Maier~(1993)~\cite{MR1225426} proved that for any constant $C > 0$ the inequality $A(n) \geq n^C$ holds on a set of positive lower density.
It is well known that $\omega(n) \sim \log \log n$ as $n \to +\infty$ over a set of natural density $1$, where $\omega(n)$ is the number of distinct prime factors of $n$ (see, e.g., \cite[Ch.~III.3]{MR3363366}).
For every $C > 1$, let $\mathcal{E}_C$ be the set of squarefree integers $n$ such that $\omega(n) \geq C \log \log n$.
Maier~(2001)~\cite{MR1830573} proved that for every $C > 2 / \log 2$ and $\varepsilon > 0$ the inequality $A(n) > \exp((\log n)^{(C \log 2) / 2 - \varepsilon})$ holds for almost all $n \in \mathcal{E}_C$.
Later, Konyagin, Maier, and Wirsing~(2004)~\cite{MR2106468} showed that, actually, such lower bound for $A(n)$ holds for all positive integers with $\omega(n) \geq C \log \log n$.
The key part of their proof is a strong upper bound on the third moment of the function $\log|\Phi_n(z)|$ over the unit circle.

\subsection{Bounds on $A(n)$ in terms of prime factors}

Felsch and Schmidt~(1968)~\cite{MR232757} and, independently, Justin~(1969)~\cite{MR241352} proved that $A(n)$ has an upper bound that does not depend on the two largest prime factors of $n$.
Let $n = p_1 \cdots p_k$, where $p_1 < \cdots < p_k$ are odd prime numbers and $k \geq 3$.
Bateman, Pomerance, and Vaughan~(1984)~\cite{MR781138} proved that $A(n) \leq M(n)$, where $M(n) := \prod_{j = 1}^{k - 2} p_j^{2^{k-j-1} - 1}$ (see also~\cite{MR2437967} for an upper bound of a similar form for $|\Phi_n(X)|$ on the unit circle).
Furthermore, they conjectured that $M(n) \leq \varphi(n)^{2^{k-1} / k - 1}$.
This conjecture was proved by Bzd{\c{e}}ga~(2012)~\cite{MR2890546}, who also proved that $A(n) \leq C_k M(n)$, where $(C_k)_{k \geq 3}$ is a sequence such that $C_k^{2^{-k}}$ converges to a constant less than $0.9541$, as $k \to +\infty$.
In the opposite direction, Bzd{\c{e}}ga~(2016)~\cite{MR3503063} proved that for every $k \geq 3$ and $\varepsilon > 0$ there exists $n$ such that $A(n) > (c_k - \varepsilon) M(n)$, where $(c_k)_{k \geq 3}$ is a sequence such that $c_k^{2^{-k}}$ converges to a constant that is about $0.71$, as $k \to +\infty$.
In~particular, this last result implies that in the upper bound on $A(n)$ the product $M(n)$ is optimal, which means that, in a precise sense, it cannot be replaced by a smaller product of $p_1, \dots, p_{k-2}$.
Furthermore, Bzd{\c{e}}ga~(2017)~\cite{MR3713089} proved several asymptotic bounds for quantities such as $A(n)$, the sum of the absolute values of the coefficients of $\Phi_n(X)$, the sum of the squares of the coefficients of $\Phi_n(X)$, and the maximum of the absolute value of $\Phi_n(X)$ on the unit circle, as $p_1 \to +\infty$ and $k$ is fixed.

\subsection{The dual function $a(j)$}

For every positive integer $j$ define
\begin{equation}\label{equ:dual-definition}
a(j) := \max_{n \,\geq\, 1} |a_n(j)| .
\end{equation}
Thus $a(j)$ is somehow a dual version of $A(n)$.
From~\eqref{equ:Phi_n-Mobius} it follows that $a_{pqn}(j) = a_n(j)$ for all prime numbers $p > q > j$ not dividing $n$.
Hence, in~\eqref{equ:dual-definition} the maximum can be replaced by a limit superior.
Erd\H{o}s and Vaughan~(1974)~\cite{MR357367} proved that $\log a(j) < 2 \tau^{1/2} j^{1/2} + C j^{3/8}$ for all positive integers $j$, where $\tau := \prod_p \left(1 - \tfrac{2}{p(p + 1)}\right)$ and $C > 0$ is a constant, and conjectured that $\log a(j) = o(j^{1/2})$ as $j \to +\infty$.
Also, they showed that $\log a(j) \gg j^{1/2} / (\log j)^{1/2}$ for all sufficiently large integers $j$.
Vaughan~(1974)~\cite{MR364141} proved that $\log a(j) \gg j^{1/2} / (\log j)^{1/4}$ for infinitely many positive integers $j$.
Montgomery and Vaughan~(1985)~\cite{MR819835} determined the order of magnitude of $\log a(j)$ by proving that that $\log a(j) \asymp j^{1/2} / (\log j)^{1/4}$ for all sufficiently large integers $j$.
Finally, Bachman~(1993)~\cite{MR1172916} proved the asymptotic formula $\log a(j) \sim C j^{1/2} / (\log j)^{1/4}$, where $C > 0$ is a constant given by a quite complicate expression.

\section{Maximum gap}

Al-Kateeb, Ambrosino, Hong, and Lee~(2021)~\cite{MR4271823} proved that $G(pm) = \varphi(m)$ for every prime number $p$ and for every squarefree positive integer $m$ with $p > m$.
This was previously numerically observed by Ambrosino, Hong, and Lee~(2017)~\cite{Ambrosino2017, AHL2017}.
The proof is based on a new divisibility property regarding a partition of $\Phi_{pm}(X)$ into ``blocks'' (see also~\cite{MR3653750, AHLblock2017}).
Furthermore, Al-Kateeb, Ambrosino, Hong, and Lee~\cite{MR4271823} conjectured that $G(pm) \leq \varphi(m)$ for every prime number $p$ and for every squarefree positive integer $m$ with $p < m$.

\section{The set of coefficients}

Suzuki~(1987)~\cite{MR931264} gave a short proof that every integer appears as the coefficient of some cyclotomic polynomial.
(Note that now this follows, for example, from Bachman's result on ternary cyclotomic polynomials with an optimally large set of coefficients~\cite{MR2053964}).
Ji and Li~(2008)~\cite{MR2459410} proved that, for each fixed prime power $p^\ell$, every integer appears as the coefficient of a cyclotomic polynomial of the form $\Phi_{p^\ell n}(X)$.
Ji, Li, and Moree~(2009)~\cite{MR2510582} generalized this result by showing that, for each fixed positive integer $m$, every integer appears as the coefficient of a cyclotomic polynomial of the form $\Phi_{mn}(X)$.
Then Fintzen~(2011)~\cite{MR2811553} determined the set $\{a_n(j) : n \equiv a \pmod d,\; j \equiv b \pmod f\}$ for any given nonnegative integers $a < d$ and $b < f$ (see also~\cite{MR3058658}).
In particular, she showed that this set is either $\mathbb{Z}$ or $\{0\}$.

Recall that $\mathcal{A}(n) := \{a_n(j) : 0 \leq j \leq \varphi(n)\}$ is the set of coefficients of $\Phi_n(X)$.
Kaplan~(2007, 2010)~\cite[Theorems~2 and 3]{MR2351667}\cite[Theorem~4]{MR2684126} proved\footnote{\cite[Theorems~2 and 3]{MR2351667} are stated with $A(n)$ in place of $\mathcal{A}(n)$, but their proofs show that the result do indeed hold with $\mathcal{A}(n)$.} the following two results regarding a kind of periodicity of $\mathcal{A}(n)$.

\begin{theorem}
Let $n$ be a binary integer, and let $p$ and $q$ be prime number greater than the largest prime factor of $n$ and such that $p \equiv \pm q \pmod n$.
Then $\mathcal{A}(pn) = \mathcal{A}(qn)$.
\end{theorem}

\begin{theorem}
Let $n$ be a positive integer, and let $p$ and $q$ be prime numbers greater than $n$ and satisfying $p \equiv q \pmod n$.
Then $\mathcal{A}(pn) = \mathcal{A}(qn)$. 
\end{theorem}

\section{Formulas for the coefficients}

Let $n > 1$ be an integer.
From~\eqref{equ:Phi_n-Mobius} it follows that
\begin{equation}\label{equ:Phi_n-infinite}
\Phi_n(X) = \prod_{d \,=\, 1}^\infty \left(1 - X^d\right)^{\mu(n/d)} ,
\end{equation}
with the convention that $\mu(x) = 0$ if $x$ is not an integer.
Therefore, each coefficient $a_n(j)$ depends only on the values $\mu(n / d)$, with $d$ a positive integer not exceeding~$j$, and using~\eqref{equ:Phi_n-infinite} one can obtain formulas for $a_n(j)$ for each fixed $j$.
For~instance, we have
\begin{align*}
a_n(1) &= -\mu(n) , \\
a_n(2) &= \tfrac1{2}\mu(n)^2 - \tfrac1{2}\mu(n) - \mu\big(\tfrac{n}{2}\big) , \\
a_n(3) &= \tfrac1{2}\mu(n)^2 - \tfrac1{2}\mu(n) + \mu(n)\mu\big(\tfrac{n}{2}\big) -\mu\big(\tfrac{n}{3}\big) .
\end{align*}
In general, M\"oller~(1970)~\cite{MR266899} proved that
\begin{equation*}
a_n(j) = \sum_{\substack{\lambda_1 + 2\lambda_2 + \,\cdots\, + j\lambda_j \,=\, j \\ \lambda_1, \dots, \lambda_j \,\geq\, 0}} \; \prod_{d \,=\, 1}^j (-1)^{\lambda_d} \binom{\mu(n / d)}{\lambda_d} ,
\end{equation*}
for every integer $j \geq 0$ (see~\cite[Lemma~4]{MR2904047} for a short proof).

Kazandzidis~(1963)~\cite{MR162793} 
and D.~H.~Lehmer~(1966)~\cite{MR197441} noted that, by Newton's identities for the symmetric elementary polynomials in terms of power sums, we have
\begin{equation*}
a_n(j) = \sum_{\substack{\lambda_1 + 2\lambda_2 + \,\cdots\, + j\lambda_j \,=\, j \\ \lambda_1, \dots,\, \lambda_j \,\geq\, 0}} \;\prod_{t \,=\, 1}^j \frac{(-c_n(t) / t)^{\lambda_t}}{\lambda_t!} ,
\end{equation*}
where
\begin{equation*}
c_n(t) := \sum_{\substack{1 \,\leq\, k \,\leq\, n \\ \gcd(n,\, k) \,=\, 1}} e^{2\pi\mathbf{i}kt / n} = \varphi(n)\,\frac{\mu(n /\!\gcd(n, t))}{\varphi(n /\!\gcd(n, t))},
\end{equation*}
is a \emph{Ramanujan's sum} and the second equality is due to H\"older~(1936)~\cite{Holder1936}.
Deaconescu and S\'{a}ndor~(1987)~\cite{DeaSan1987} (see also~\cite[pp.~258--259]{MR2119686}) gave another formula for $a_n(j)$ in terms of a determinant involving Ramanujan's sums.
Furthermore, Eaton~(1939)~\cite{MR1563935} proved a formula for $a_n(j)$ in terms of a sum having each addend either equal to $-1$ or $+1$ depending on a quite involved rule.

Grytczuk and Tropak~(1991)~\cite{MR1151850} provided another method to compute $a_n(j)$, which makes use of the recurrence
\begin{equation*}
a_n(j) = -\frac{\mu(n)}{j}\sum_{i \,=\, 0}^{j - 1} a_n(i) \, \mu(\gcd(n, j - i)) \, \varphi(\gcd(n, j - i)) , \quad \text{ for } j > 0,
\end{equation*}
with $a_n(0) = 1$.
By using this method, they found for $m=\pm 2, \dots, \pm 9$, and $10$ the minimal positive integer $j$ for which there exists a positive integer $n$ such that $a_n(j) = m$.

Herrera-Poyatos and Moree~(2021)~\cite{HerrMoree2021} wrote a survey on formulas for $a_n(j)$ involving Bernoulli numbers, Stirling numbers, and Ramanujan's sums.
Also, they introduced a new uniform approach that makes possible to provide shorter proofs for some of such formulas and also to derive new ones.

\section{Miscellaneous results}
Carlitz~(1967)~\cite{MR217044} proved some asymptotic formulas involving the sum of squares of the coefficients of $\Phi_n(X)$.
Endo~(1974)~\cite{MR396399} proved that $7$ is the minimal nonnegative integer $j$ such that $|a_n(j)| > 1$ for some positive integer $n$.
Dresden~(2004)~\cite{MR3952408} proved that for every $n \geq 3$ the middle coefficient of $\Phi_n(X)$ is either $0$, and in such a case $n$ is a power of $2$, or an odd integer.
Dunand~(2012)~\cite{MR2911275} studied the coefficients of the inverse of $\Phi_m(X)$ modulo $\Phi_n(X)$, where $m$ and $n$ are distinct divisors of $pq$, with $p < q$ primes, and discussed an application to torus-based cryptography.
Musiker and Reiner~(2014)~\cite{MR3176609} gave an interpretation of $a_n(j)$ as the torsion order in the homology of certain simplicial complexes.
An alternative proof of this results was given by Meshulam~(2012)~\cite{MR2902701}.
Chu~(2021)~\cite{MR4217760} gave necessary conditions on $n$ so that the powers of positive, respectively negative, coefficients of $\Phi_n(X)$ are in arithmetic progression.
For all integers $j, v \geq 0$, let 
\begin{equation*}
\overline{a}(j) := \lim_{N \to +\infty} \frac1{N} \sum_{n \,\leq\, N} a_n(j)
\end{equation*}
be the average value of the $j$th coefficient of the cyclotomic polynomials, and let
\begin{equation*}
\delta(j, v) := \lim_{N \to +\infty} \frac1{N} \#\{n \leq N : a_n(j) = v\} ,
\end{equation*}
be the frequency that such coefficient is equal to $v$.
M\"{o}ller~(1970)~\cite{MR266899} proved that $\overline{a}(j) = \frac{6}{\pi^2} e_j$ for every integer $j \geq 1$, where $e_j > 0$ is a rational number.
Gallot, Moree, and Hommersom~(2011)~\cite{MR2904047} derived explicit formulas for $\overline{a}(j)$ and $\delta(j, v)$.
Also, they verified that $f_j := e_j j \prod_{p \leq j}(p + 1)$ is an integer for every positive integer $j \leq 100$, and asked whether it is true in general.
Gong~(2009)~\cite{MR2560843} proved that indeed every $f_j$ is an integer, and also showed that, for every integer $m$, we have that $m \mid f_j$ for every sufficiently large $j$.

\section{Algorithms and numerical data}

Arnold and Monagan~(2011)~\cite{MR2813365} presented three algorithms for computing the coefficients of the $n$th cyclotomic polynomial, and wrote a fast implementation using machine-precision arithmetic.
The first algorithm computes $\Phi_n(X)$ by a series of polynomial divisions using Lemma~\ref{lem:basic}'s \ref{lem:basic:2}.
This method is well known~\cite{MR1344932}, but Arnold and Monagan optimized the polynomial division by way of the discrete Fast Fourier Transform.
The second algorithm computes $\Phi_n(X)$ as a quotient of sparse power series using~\eqref{equ:Phi_n-Mobius}.
In such algorithm, $\Phi_n(X)$ is treated as a truncated power series.
Multiplication of a truncated power series by $X^d - 1$ is easy, and division by $X^d - 1$ is equivalent to multiplication by the power series $-\sum_{j=0}^\infty X^{dj}$.
This algorithm was further improved in a subsequent work~\cite{Arnold2010}.
The third algorithm, which they called the ``big prime algorithm'', generates the terms of $\Phi_n(X)$ sequentially, in a manner which reduces the memory cost.

\begin{figure}[h]
 \centering
 \includegraphics[width=\textwidth]{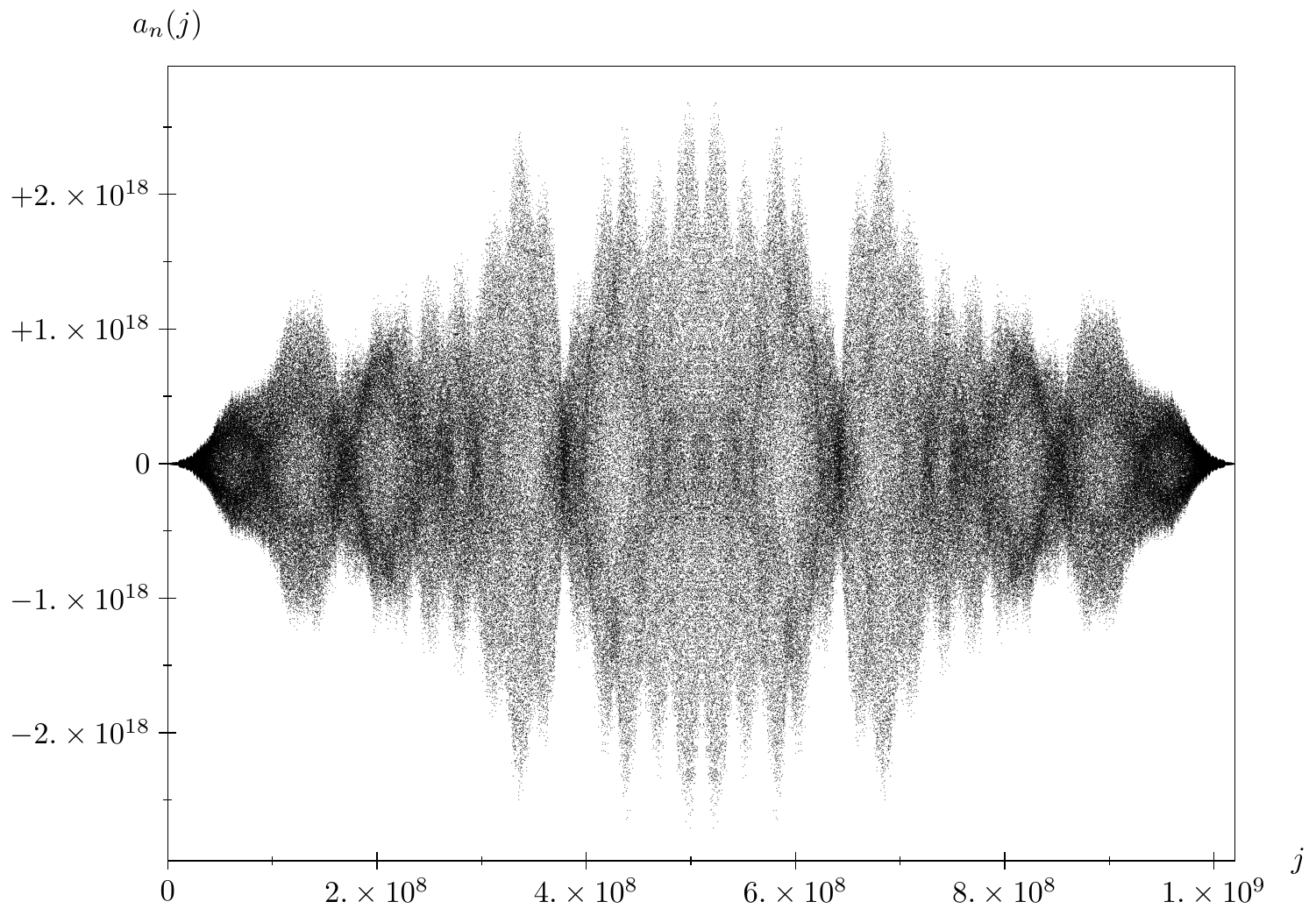}
 \caption{A plot of the coefficients of $\Phi_n(X)$ for $n = 3\cdot 5 \cdot 7 \cdot 11 \cdot 13 \cdot 17 \cdot 19 \cdot 23 \cdot 29$.
 The $\varphi(n) + 1 = 1,021,870,081$ coefficients were computed using the program 
 \texttt{SPS4\_64} of Arnold and Monagan~\cite{AMdata}.
 Then the plot was produced by selecting a random sample of $500,000$ coefficients.
 }
 \label{fig:Phi_3-29}
\end{figure}

With their implementation, Arnold and Monagan produced a large amount of data on the coefficients of $\Phi_n(X)$ for $n$ in the range of billions~\cite{AMdata}.
For instance, they found the minimal positive integer $n$ such that $A(n)$ is greater than $n$, $n^2$, $n^3$, and $n^4$, respectively.
Also, they computed $A(n)$ when $n$ is equal to the product of the first $9$ odd prime numbers. 
(Partial computations on the cases of $n$ equal to the product of the first $7$ and $8$ odd prime numbers were previously done by Koshiba~(1998, 2000)~\cite{MR1718402, MR1833785}.)

Other numerical data on the coefficients of the cyclotomic polynomial can be found on the Online Encyclopedia of Integer Sequences~\cite{OEIS}.
See for instance sequences \texttt{A117223}, \texttt{A117318}, \texttt{A138474}, and \texttt{A138475} of Noe.

\section{Relatives of cyclotomic polynomials}

In this section we collect results regarding the coefficients of polynomials that are closely related to cyclotomic polynomials.

\subsection{Inverse cyclotomic polynomials}

Let $n$ be a positive integer.
The \emph{$n$th inverse cyclotomic polynomial} $\Psi_n(X)$ is defined as the monic polynomial whose roots are exactly the nonprimitive $n$th roots of unity, that~is,
\begin{equation}\label{equ:Psi_n-definition}
\Psi_n(X) := \prod_{\substack{1 \,\leq\, k \,\leq\, n \\ \gcd(n,\, k) \,>\, 1}} \left(X - e^{2\pi\mathbf{i}k / n}\right) = \frac{X^n - 1}{\Phi_n(X)} .
\end{equation}
Note that $\Phi_n(X)$ has degree $n - \varphi(n)$.
From~\eqref{equ:Xn-1} and~\eqref{equ:Psi_n-definition} it follows that
\begin{equation*}
\Psi_n(X) = \prod_{\substack{d \,\mid\, n \\ d \,>\, 1}} \Phi_d(X) .
\end{equation*}
In~particular, $\Psi_n(X)$ has integer coefficients.
Moreover, from~\eqref{equ:Psi_n-definition} we get that
\begin{equation*}
\frac1{\Phi_n(X)} = \frac{\Psi_n(X)}{X^n - 1} = -\Psi_n(X) \sum_{j \,=\, 0}^\infty X^{nj} .
\end{equation*}
Thus, the Taylor coefficients of $1 / \Phi_n(X)$ are purely periodic, and the period consists of the $n - \varphi(n) + 1$ coefficients of $-\Phi_n(X)$ followed by $\varphi(n) - 1$ zeros.
The next lemma collects some basic identities, which follows easily from Lemma~\ref{lem:basic} and~\eqref{equ:Psi_n-definition}.

\begin{lemma}\label{lem:Psi_n-basic}
For every positive integer $n$ and every prime number $p$, we have that:
\begin{enumerate}
\item $\Psi_{pn}(X) = \Psi_n(X^p)$ if $p \mid n$;
\item $\Psi_{pn}(X) = \Phi_n(X) \Psi_n(X^p)$ if $p \nmid n$;
\item\label{lem:Psi_n-basic3} $\Psi_{2n}(X) = (-1)^{\varphi(n)} (1 - X^n) \Psi_n(-X)$ if $2 \nmid n$;
\item\label{lem:Psi_n-basic4} $\Psi_n(X) = \Psi_{\rad(n)}(X^{n / \!\rad(n)})$;
\item $\Psi_n(1/X) = -X^{-(n - \varphi(n))} \Psi_n(X)$ if $n > 1$.
\end{enumerate}
\end{lemma}

Similarly to cyclotomic polynomials, in light of Lemma~\ref{lem:Psi_n-basic}'s \ref{lem:Psi_n-basic3} and \ref{lem:Psi_n-basic4}, for the purpose of studying the coefficients of the inverse cyclotomic polynomial $\Psi_n(X)$ it suffices to consider only odd squarefree integers $n$.
For a squarefree positive integer $n$, the inverse cyclotomic polynomial $\Psi_n(X)$ is \emph{binary}, \emph{ternary}, \dots if the number of prime factors of $n$ is $2$, $3$, \dots.
The \emph{order} of $\Psi_n(X)$ is the number of prime factors of $n$.
It is easy to check that $\Psi_1(X) = 1$, $\Psi_p(X) = X - 1$, and
\begin{equation*}
\Psi_{pq}(X) = X^{p + q - 1} + X^{p + q - 2} + \cdots + X^q - X^{p - 1} - X^{p - 2} - \cdots - 1
\end{equation*}
for all prime numbers $p < q$.
Hence, the simplest nontrivial case in the study of the coefficients of $\Psi_n(X)$ occurs for ternary $n$.

Let $C(n)$ denote the height of $\Psi_n(X)$.
Moree~(2009)~\cite{MR2488596} proved that
\begin{equation*}
C(pqr) \leq \left\lfloor \frac{(p - 1)(q - 1)}{r}\right\rfloor + 1 \leq p - 1 ,
\end{equation*}
for all odd primes $p < q < r$.
Also, he showed that $C(pqr) = p - 1$ if and only if $q \equiv r \equiv \pm 1 \pmod p$ and $r < \tfrac{p - 1}{p - 2}(q - 1)$.
Furthermore, he provided several results on flat inverse cyclotomic polynomials.
For instance, he showed that $\Psi_{15r}(X)$ and $\Psi_{21r}(X)$ are flat, for every prime $p$, and that $\Psi_{pqr}(X)$ is flat for all primes $p < q$ and $r > (p - 1)(q - 1)$.
Furthermore, he proved that every integer appears as the coefficient of some inverse cyclotomic polynomial.
Bzd{\c{e}}ga~(2014)~\cite{MR3195385} proved a formula for $C(pqr)$ in the case in which $r = \alpha p + \beta q \leq \varphi(pq)$ for some positive integers $\alpha, \beta$.
Using such formula, he gave necessary and sufficient conditions for $\Psi_{pqr}(X)$ being flat in such a case.
Hong, Lee, Lee, and Park~(2012)~\cite{MR2944756} proved that $G(\Psi_{pqr}) = 2qr - \deg(\Psi_{pqr})$ for all odd primes $p < q < r$ such that $q > 4(p - 1)$ or $r > p^2$.
Also, they gave lower and upper bound for $G(\Psi_{pqr})$ for general $\Psi_{pqr}$.
In general, many papers regarding the coefficients of cyclotomic polynomials also provide related results for the coefficients of inverse cyclotomic polynomials~\cite{Arnold2010, MR2813365, MR2890546, MR3459568, MR2811553, MR2568054, MR2904047, HerrMoree2021, MR3922607}.

\subsection{Divisors of \texorpdfstring{$X^n - 1$}{X\textnsuperior{} -- 1}}

A natural generalization of the study of the coefficients of $\Phi_n(X)$ is the study of the coefficients of divisors of $X^n - 1$.
Note that, in light of~\eqref{equ:Xn-1} and the irreducibility of cyclotomic polynomials, $X^n - 1$ has $2^{\tau(n)}$ monic divisors in $\mathbb{Z}[X]$, where $\tau(n)$ is the number of (positive) divisors of $n$, which are given by products of distinct cyclotomic polynomials $\Phi_d(X)$ with $d$ a divisor of $n$.
Let $B(n)$ be the maximum height of the monic divisors of $X^n - 1$.
Justin~(1969)~\cite{MR241352} showed that $B(n)$ has an upper bound that is independent from the largest prime factor of $n$.
Pomerance and Ryan~(2007)~\cite{MR2342677} proved that
\begin{equation*}
\limsup_{n \to +\infty} \frac{\log \log B(n)}{\log n / \log \log n} = \log 3 .
\end{equation*}
Furthermore, they showed that $B(pq) = p$ for all primes $p < q$, and that $B(n) = 1$ if and only if $n$ is a prime power.
Kaplan~(2009)~\cite{MR2549523} proved that $B(p^2 q) = \min\{p^2, q\}$ for all distinct primes $p$ and $q$, and that
\begin{equation*}
\tfrac1{3}(3p^2 q - p^3 + 7p - 6) \leq B(pqr) \leq p^2 q^2 ,
\end{equation*}
for all primes $p < q < r$.
Moreover, letting $n = p_1^{e_1} \cdots p_k^{e_k}$, where $p_1 < \cdots < p_k$ are prime numbers, $e_1, \dots, e_k$ are positive integers, and $k \geq 2$.
Kaplan proved the upper bound $B(n) < \prod_{j=1}^{k-1} p_i^{4\cdot 3^{k-2} E - e_j}$, where $E := \prod_{j=1}^k e_j$.
Bzd{\c{e}}ga~(2012)~\cite{MR2890546} showed that $B(n) < (C + o(1))^{3^k} n^{(3^k - 1)/(2k) - 1}$, as $k \to +\infty$, where $C < 0.9541$ is an effectively computable constant.
Zhang~(2019)~\cite{MR4007431} improved Kaplan's bound to $B(n) < \left(\tfrac{2}{5}\right)^{\prod_{j=2}^k e_j} p_i^{4\cdot 3^{k-2} E - e_j}$.
Ryan, Ward, and Ward~(2010)~\cite{MR2763271} proved that $B(n) \geq \min\{u, v\}$ whenever $n = uv$, where $u$ and $v$ are coprime positive integers.
In~particular, this implies that $B(n) \geq \min\{p_1^{e_1}, \cdots, p_k^{e_k}\}$.
Furthermore, they made several conjectures on $B(n)$, for $n$ having two, three, or four prime factors, based on extensive numerical computations.
Some of these conjectures were proved by Wang~(2015)~\cite{MR3280943}.
In~particular, he showed that for all odd primes $p < q < r$ and every positive integer $b$ we have that: $B(pq^b)$ is divisible by $p$, $B(2q^b) = 2$, if $b \geq 3$ then $B(pq^b) > p$, and if $q \equiv \pm r \pmod p$ and $b \leq 5$ then $B(pq^b) = B(pr^b)$.
Thompson~(2011)~\cite{MR2988080} proved that $B(n) \leq n^{\tau(n) f(n)}$ for almost all positive integers $n$, where $f(n)$ is any function such that $f(n) \to +\infty$ as $n \to +\infty$.
Decker and Moree~(2013)~\cite{MR3071668} (see also the extended version~\cite{MR3071668_Extended}) determined the set of coefficients of each of the $64$ divisors of $X^{p^2 q} - 1$, where $p$ and $q$ are distinct primes. 
In~particular, their result shows that for most of the divisors the set of coefficients consists of consecutive integers. 
Moreover, they proved that if $f_e$ is the number of flat divisors of $X^{p^e q} - 1$, for each integer $e \geq 1$, then $f_{e + 1} \geq 2f_e + 2^{e + 2} - 1$.

For each integer $r \geq 1$, let $B(r, n)$ be the maximum of the absolute value of the coefficient of $X^r$ in $f(X)$, as $f(X)$ ranges over the monic divisors of $X^n - 1$.
Somu~(2016)~\cite{MR3504047} gave upper and lower bounds for $B(r, n)$ that imply
\begin{equation*}
\limsup_{n \to +\infty} \frac{\log B(r, n)}{\log n / \log \log n} = r \log 2 .
\end{equation*}
In the same work, Somu proved that if $\ell$ and $m$ are positive integers, then there exist a positive integer $n$ and a monic divisor $f(X)$ of $X^n - 1$ having exactly $m$ irreducible factors such that each integers in $[-\ell, \ell]$ appears among the coefficients of $f(X)$.
Moreover, he showed that for all integers $c_1, \dots, c_r$ there exist a positive integer $n$ and a divisor $f(X) = \sum_{j = 1}^{\deg(f)} d_j X^j$, with $f_i \in \mathbb{Z}$, of $X^n - 1$ such that $d_i = c_i$ for $i=1,\dots, r$.
Later Somu~(2017)~\cite{MR3541693} proved that the set of such $n$ has positive natural density.

\subsection{Inclusion-exclusion polynomials}\label{sec:inclusion-exclusion}

Inclusion-exclusion polynomials were introduced by Bachman~(2010)~\cite{MR2798626} as a kind of combinatorial generalization of cyclotomic polynomials.
Let bold letters $\bm{n}, \bm{d}, \dots$ denote finite sets of pairwise coprime integers greater than $1$.
Furthermore, for each $\bm{n} = \{n_1, \dots, n_k\}$, where $n_1, \dots, n_k > 1$ are pairwise coprime integers, put $\|\bm{n}\| := n_1 \cdots n_k$, $\mu(\bm{n}) := (-1)^k$, and $\varphi(\bm{n}) := \prod_{i=1}^k (n_i - 1)$.
The \emph{$\bm{n}$th inclusion-exclusion polinomial} is defined as
\begin{equation}\label{equ:inclusion-exclusion-definition}
\Phi_{\bm{n}}(X) = \prod_{\bm{d} \,\subseteq\, \bm{n}} \left(X^{\|\bm{n}\| / \|\bm{d}\|} - 1\right)^{\mu(\bm{d})} .
\end{equation}
Note the striking resemblance of~\eqref{equ:Phi_n-Mobius} and~\eqref{equ:inclusion-exclusion-definition}.
In~particular, we have that 
\begin{equation*}
\Phi_{\{p_1, \dots, p_k\}}(X) = \Phi_{p_1 \cdots p_k}(X) ,
\end{equation*}
for all prime numbers $p_1 < \cdots < p_k$.

Many results regarding cyclotomic polynomials can be generalized to inclusion-exclusion polynomials, and it might be even more natural to prove them directly for inclusion-exclusion polynomials~\cite{MR3206393, MR3459568, MR4217760, MR3208875, MR3295662}.
Also, the \emph{$\bm{n}$th inverse inclusion-exclusion polynomial}, defined by $\Psi_{\bm{n}}(X) := (X^{\|\bm{n}\|} - 1) / \Phi_{\bm{n}}(X)$, has been studied~\cite{MR3195385}.

The following theorem summarizes the basic properties of inclusion-exclusion polynomials, including the fact that they are indeed polynomials~\cite{MR2798626}.

\begin{theorem}
For every $\bm{n} = \{n_1, \dots, n_k\}$, where $n_1, \dots, n_k > 1$ are pairwise coprime integers, we have that
\begin{equation*}
\Phi_{\bm{n}}(X) = \prod_{\omega} (X - \omega) ,
\end{equation*}
where $\omega$ runs over the $\|\bm{n}\|$th roots of unity satisfying $\omega^{\|\bm{n}\| / n_i} \neq 1$ for all $i=1,\dots,k$. 
Moreover, the degree of $\Phi_{\bm{n}}(X)$ is equal to $\varphi(\bm{n})$ and it holds
\begin{equation*}
\Phi_{\bm{n}}(X) = \prod_{d} \Phi_d(X) ,
\end{equation*}
where $d$ runs over the divisors of $\|\bm{n}\|$ such that $(d, n_i) > 1$ for every $i=1,\dots,k$.
In~particular, $\Phi_{\bm{n}}(X)$ has integer coefficients.
\end{theorem}

Let $p, q, r, s$ be pairwise coprime integers greater than $1$.
Bachman~(2010)~\cite{MR2798626} proved that the set of coefficients of every \emph{ternary inclusion-exclusion polynomial} $\Phi_{\{p,q,r\}}(X)$ consists of consecutive integers and, assuming $p < q < r$, it depends only on the residue class of $r$ modulo $pq$.
Let $A(p, q, r)$ denote the height of $\Phi_{\{p, q, r\}}(X)$.
Bachman and Moree~(2011)~\cite{MR2798664} showed that, if $r \equiv \pm s \pmod {pq}$ and $r > \max\{p,q\} > s \geq 1$, then
\begin{equation*}
A(p,q,s) \leq A(p,q,r) \leq A(p,q,s) + 1 .
\end{equation*}
For every $\bm{n} = \{n_1, \dots, n_k\}$, where $n_1 < \cdots < n_k$ are pairwise coprime integers greater than $1$, let $A(\bm{n})$ be the height of $\Phi_{\bm{n}}(X)$ and put $M(\bm{n}) := \prod_{j = 1}^{k-2} n_j^{2^{k-j-1}-1}$.
Also, let $D_k$ be the smallest real number for which the inequality $A(\bm{n}) \leq D_k M(\bm{n})$ holds for all sufficiently large $n_1$.
Bzd{\c{e}}ga~(2013)~\cite{MR3141834} proved that $(C_1 + o(1))^{2^k} < D_k < (C_2 + o(1))^{2^k}$, as $k \to \infty$, where $C_1, C_2 > 0$ are constants, with $C_1 \approx 0.5496$ and $C_2 \approx 0.9541$.
Furthermore, Liu~(2014)~\cite{MR3268704} studied the polynomial obtained by restricting~\eqref{equ:inclusion-exclusion-definition} to the sets $\bm{d}$ with at most two elements.

\subsection{Unitary cyclotomic polynomials}

Let $n$ be a positive integer.
A \emph{unitary divisor} of $n$ is a divisor $d$ of $n$ such that $d$ and $n/d$ are relatively prime.
Moree and T\'{o}th~(2020)~\cite{MR4142502} defined the \emph{$n$th unitary cyclotomic polynomial} as
\begin{equation*}
\Phi_n^*(X) := \prod_{\substack{1 \,\leq\, k \,\leq\, n \\[2pt] \gcd^*\!(n,\, k) \,=\, 1}} \left(X - e^{2\pi\mathbf{i}k / n}\right) ,
\end{equation*}
where $\gcd^*(n, k)$ denotes the maximum unitary divisor of $n$ which is a divisor of $k$.
\mbox{It can be proved} that $\Phi_n^*(X)$ has integer coefficients.
Moreover, the following analogs of~\eqref{equ:Xn-1} and~\eqref{equ:Phi_n-Mobius} holds:
\begin{equation*}
X^n - 1 = \prod_{d \,\mid\mid\, n} \Phi_d^*(X) ,
\end{equation*}
where $d \mid\mid n$ means that $d$ is a unitary divisor of $n$, and
\begin{equation*}
\Phi_n^*(X) = \prod_{d \,\mid\mid\, n} \left(X^{n / d} - 1\right)^{\mu^*(d)} ,
\end{equation*}
where $\mu^*(n) := (-1)^{\omega(n)}$.
Every unitary cyclotomic polynomial can be written as an inclusion-exclusion polynomial, precisely $\Phi_n^*(X) = \Phi_{\{p_1^{e_1}, \dots, \, p_k^{e_k}\}}(X)$ for $n = p_1^{e_1} \cdots p_k^{e_k}$, where $p_1 < \cdots < p_k$ are prime numbers and $e_1, \dots, e_k$ are positive integers.
Furthermore, every unitary cyclotomic polynomial is equal to a product of cyclotomic polynomials:
\begin{equation*}
\Phi_n^*(X) = \prod_{\substack{d \,\mid\, n \\ \rad(d) \,=\, \rad(n)}} \Phi_d(X) .
\end{equation*}
These and other properties of unitary cyclotomic polynomials were proved by Moree and T\'{o}th~\cite{MR4142502}.
Jones, Kester, Martirosyan, Moree, T\'{o}th, White, and Zhang~(2020)~\cite{MR4121365} proved that, given any positive integer $m$, every integer appears as a coefficient of $\Phi_{mn}^*(X)$, for some positive integer $n$.
Also, they showed the analog result for coefficients of the \emph{inverse unitary cyclotomic polynomial} $\Psi_n^*(X) := (X^n - 1) / \Phi_n^*(X)$.
Bachman~(2021)~\cite{Bachman2021} proved that, fixed three distinct odd primes $p, q, r$ and $\varepsilon > 0$, for every sufficiently large positive integer $a$, depending only on $\varepsilon$, there exist positive integers $b$ and $c$ such that the the set of coefficients of $\Phi_{p^a q^b r^c}^*(X)$ contains all the integers in the interval $\big[{-(\tfrac1{4}-\varepsilon)p^a}, (\tfrac1{4}-\varepsilon)p^a\big]$.
As a consequence, every integer appears as the coefficient of some ternary unitary cyclotomic polynomial.
Furthermore, he provided an infinite family of ternary unitary cyclotomic polynomials $\Phi_{p^a q^b r^c}^*(X)$ whose sets of coefficients consist of all the integers in $[-(p^a - 1)/2, (p^a + 1)/2]$, and he pointed out that this interval is as large as possible.

\newpage
\cleardoublepage
\phantomsection
\addcontentsline{toc}{section}{References}
\bibliographystyle{amsplain}

\end{document}